\documentclass [10pt,]{article}
\usepackage{amssymb,amsmath}
\usepackage{indentfirst}
\numberwithin{equation}{section}


\def\BigRoman{\uppercase\expandafter{\romannumeral\number\count 255 }}
\def\Romannumeral{\afterassignment\BigRoman\count255=}

\newtheorem{Theorem}{Theorem}
\newtheorem{Lemma}{Lemma}[section]

\newtheorem{Corollary}{Corollary}

\begin{document}%

\title{A neighborhood condition for fractional ID-$[a,b]$-factor-critical graphs
}
\author{Sizhong Zhou  \thanks{Corresponding
 author. E-mail address: zsz\_cumt@163.com(S. Zhou)}\\
School of Mathematics and Physics\\
Jiangsu University of Science and Technology\\
Mengxi Road 2, Zhenjiang, Jiangsu 212003, P. R. China\\
Fan Yang\\
School of Mathematics and Physics\\
Jiangsu University of Science and Technology\\
Mengxi Road 2, Zhenjiang, Jiangsu 212003, P. R. China\\
Zhiren Sun\\
School of Mathematical Sciences, Nanjing Normal University\\
Nanjing, Jiangsu 210046, P. R. China\\
 }
\date{}
\maketitle

\begin{abstract}Let $G$ be a graph of order $n$, and let $a$
and $b$ be two integers with $1\leq a\leq b$. Let $h:
E(G)\rightarrow [0,1]$ be a function. If $a\leq\sum_{e\ni
x}h(e)\leq b$ holds for any $x\in V(G)$, then we call $G[F_h]$ a
fractional $[a,b]$-factor of $G$ with indicator function $h$ where
$F_h=\{e\in E(G): h(e)>0\}$. A graph $G$ is fractional
independent-set-deletable $[a,b]$-factor-critical (in short,
fractional ID-$[a,b]$-factor-critical) if $G-I$ has a fractional
$[a,b]$-factor for every independent set $I$ of $G$. In this
paper, it is proved that if $n\geq\frac{(a+2b)(2a+2b-3)+1}{b}$,
$\delta(G)\geq\frac{bn}{a+2b}+a$ and $|N_G(x)\cup
N_G(y)|\geq\frac{(a+b)n}{a+2b}$ for any two nonadjacent vertices
$x,y\in V(G)$, then $G$ is fractional ID-$[a,b]$-factor-critical.
Furthermore, it is shown that this result is best possible in some
sense.\\


 {\bf (2010) Mathematics Subject Classification:} 05C70, 05C72, 05C35
\end{abstract}

\section{Introduction}
The graphs considered here will be finite undirected graphs
without loops or multiple edges. Let $G$ be a graph. We denote by
$V(G)$ and $E(G)$ the set of vertices and the set of edges of $G$,
respectively. For any $x\in V(G)$, we denote the degree of $x$ in
$G$ by $d_G(x)$. We write $N_G(x)$ for the set of vertices
adjacent to $x$ in $G$, and $N_G[x]$ for $N_G(x)\cup\{x\}$. For
$S\subseteq V(G)$, we use $G[S]$ to denote the subgraph of $G$
induced by $S$, and $G-S=G[V(G)\setminus S]$. Let $S$ and $T$ be
two disjoint vertex subsets of $G$, we denote the number of edges
from $S$ to $T$ by $e_G(S,T)$. We denote by $\delta(G)$ the
minimum degree of $G$. For any nonempty subset $S$ of $V(G)$, let
$$
N_G(S)=\bigcup_{x\in S}N_G(x).
$$
If $G$ and $H$ are disjoint graphs, the join and the union are
denoted by $G\vee H$ and $G\cup H$, respectively.

A factor of a graph $G$ is a spanning subgraph of $G$. Let $a$ and
$b$ be two positive integers with $1\leq a\leq b$. Then a factor
$F$ is an $[a,b]$-factor if $a\leq d_F(x)\leq b$ for each $x\in
V(G)$. If $a=b=k$, then an $[a,b]$-factor is called a $k$-factor.
If $k=1$, then a $1$-factor is also called a perfect matching. A
graph $G$ is factor-critical [1] if $G-x$ has a perfect matching
for each $x\in V(G)$. The concept of the factor-critical graph was
generalized to the ID-factor-critical graph [2]. We say that $G$
is independent-set-deletable factor-critical (shortly,
ID-factor-critical) if for every independent set $I$ of $G$ which
has the same parity with $|V(G)|$, $G-I$ has a perfect matching.
Obviously, every ID-factor-critical graph with odd vertices is
factor-critical.

Let $h: E(G)\rightarrow [0,1]$ be a function. If $a\leq\sum_{e\ni
x}h(e)\leq b$ holds for any $x\in V(G)$, then we call $G[F_h]$ a
fractional $[a,b]$-factor of $G$ with indicator function $h$ where
$F_h=\{e\in E(G): h(e)>0\}$. If $a=b=k$, then a fractional
$[a,b]$-factor is a fractional $k$-factor. A fractional 1-factor
is also called a fractional perfect matching. A graph $G$ is
fractional independent-set-deletable $[a,b]$-factor-critical (in
short, fractional ID-$[a,b]$-factor-critical) [3] if $G-I$ has a
fractional $[a,b]$-factor for every independent set $I$ of $G$. If
$a=b=k$, then a fractional ID-$[a,b]$-factor-critical graph is a
fractional ID-$k$-factor-critical graph. If $k=1$, then a
fractional ID-$k$-factor-critical graph is called a fractional
ID-factor-critical graph.

Many authors have investigated graph factors [4-9]. Chang, Liu and
Zhu [10] showed a minimum degree condition for a graph to be a
fractional ID-$k$-factor-critical graph. Zhou, Xu and Sun [11]
obtained an independence number and minimum degree condition for
graphs to be fractional ID-$k$-factor-critical graphs. Zhou, Sun
and Liu [3] obtained a minimum degree condition for a graph to be
a fractional ID-$[a,b]$-factor-critical graph. In this paper, we
proceed to study fractional ID-$[a,b]$-factor-critical graphs, and
obtain a neighborhood condition for a graph to be fractional
ID-$[a,b]$-factor-critical. The main result is the following
theorem.
\begin{Theorem} Let  $1\leq a\leq b$ be two
integers, and let $G$ be a graph of order $n$ with
$n\geq\frac{(a+2b)(2a+2b-3)+1}{b}$, and
$\delta(G)\geq\frac{bn}{a+2b}+a$. If $|N_G(x)\cup
N_G(y)|\geq\frac{(a+b)n}{a+2b}$ for any two nonadjacent vertices
$x,y\in V(G)$, then $G$ is fractional ID-$[a,b]$-factor-critical.
\end{Theorem}

If $a=b=k$ in Theorem 1, then we obtain the following corollary.
\begin{Corollary} Let  $k\geq1$ be an
integer, and let $G$ be a graph of order $n$ with $n\geq12k-8$,
and $\delta(G)\geq\frac{n}{3}+k$. If $|N_G(x)\cup
N_G(y)|\geq\frac{2n}{3}$ for any two nonadjacent vertices $x,y\in
V(G)$, then $G$ is fractional ID-$k$-factor-critical.
\end{Corollary}

If $k=1$ in Corollary 1, then we get the following corollary.
\begin{Corollary} Let $G$ be a graph of order $n$ with $n\geq4$,
and $\delta(G)\geq\frac{n}{3}+1$. If $|N_G(x)\cup
N_G(y)|\geq\frac{2n}{3}$ for any two nonadjacent vertices $x,y\in
V(G)$, then $G$ is fractional ID-factor-critical.
\end{Corollary}

\section{The Proof of Theorem 1}
In order to prove Theorem 1, we rely heavily on the following
lemma.
\begin{Lemma}$^{[12]}$ Let $G$ be a graph. Then $G$ has a
fractional $[a,b]$-factor if and only if for every subset $S$ of
$V(G)$,
$$
\delta_G(S,T)=b|S|+d_{G-S}(T)-a|T|\geq0,
$$
where $T=\{x:x\in V(G)\setminus S, d_{G-S}(x)\leq a\}$ and
$d_{G-S}(T)=\sum_{x\in T}d_{G-S}(x)$.
\end{Lemma}

{\bf Proof of Theorem 1.} \ Let $X$ be an independent set of $G$
and $H=G-X$. In order to complete the proof of Theorem 1, we need
only to prove that $H$ has a fractional $[a,b]$-factor. By
contradiction, we suppose that $H$ has no fractional
$[a,b]$-factor. Then by Lemma 2.1, there exists some subset
$S\subseteq V(H)$ such that
$$
\delta_H(S,T)=b|S|+d_{H-S}(T)-a|T|\leq-1,\eqno(1)
$$
where $T=\{x:x\in V(H)\setminus S, d_{H-S}(x)\leq a\}$. We first
prove the following claims.

{\bf Claim 1.} \ \ $|X|\leq\frac{bn}{a+2b}$.

{\bf Proof.} \ Since $n\geq\frac{(a+2b)(2a+2b-3)+1}{b}$, the
inequality holds for $|X|=1$. In the following we may assume
$|X|\geq2$. In terms of the condition of Theorem 1, there exist
$x,y\in X$ such that $|N_G(x)\cup N_G(y)|\geq\frac{(a+b)n}{a+2b}$.
Since $X$ is independent, we obtain $X\cap(N_G(x)\cup
N_G(y))=\emptyset$. Thus, we have
$$
|X|+\frac{(a+b)n}{a+2b}\leq|X|+|N_G(x)\cup N_G(y)|\leq n,
$$
which implies
$$
|X|\leq n-\frac{(a+b)n}{a+2b}=\frac{bn}{a+2b}.
$$
This completes the proof of Claim 1.

{\bf Claim 2.} \ \ $\delta(H)\geq a$.

{\bf Proof.} \ Note that $H=G-X$. Combining this with Claim 1, we
obtain
$$
\delta(H)\geq\delta(G)-|X|\geq(\frac{bn}{a+2b}+a)-\frac{bn}{a+2b}=a.
$$
The proof of Claim 2 is complete.

{\bf Claim 3.} \ \ $|T|\geq b+1$.

{\bf Proof.} \ If $|T|\leq b$, then from Claim 2 and
$|S|+d_{H-S}(x)\geq d_H(x)\geq\delta(H)$ for each $x\in T$, we
have
\begin{eqnarray*}
\delta_H(S,T)&=&b|S|+d_{H-S}(T)-a|T|\geq|T||S|+d_{H-S}(T)-a|T|\\
&=&\sum_{x\in T}(|S|+d_{H-S}(x)-a)\geq\sum_{x\in
T}(\delta(H)-a)\geq0,
\end{eqnarray*}
which contradicts (1). This completes the proof of Claim 3.

{\bf Claim 4.} \ \ $a|T|>b|S|$.

{\bf Proof.} \ If $a|T|\leq b|S|$, then from (1) we obtain
$$
-1\geq\delta_H(S,T)=b|S|+d_{H-S}(T)-a|T|\geq b|S|-a|T|\geq0,
$$
it is a contradiction. This completes the proof of Claim 4.

{\bf Claim 5.} \ \ $|S|+|X|<\frac{(a+b)n}{a+2b}$.

{\bf Proof.} \ According to Claim 1, Claim 4 and $|S|+|T|+|X|\leq
n$, we have
$$
an\geq
a|S|+a|T|+a|X|>a|S|+b|S|+a|X|=(a+b)(|S|+|X|)-b|X|\geq(a+b)(|S|+|X|)-\frac{b^{2}n}{a+2b},
$$
which implies
$$
|S|+|X|<\frac{(a+b)n}{a+2b}.
$$
The proof of Claim 5 is complete.

In view of Claim 3, $T\neq\emptyset$. Define
$$
h_1=\min\{d_{H-S}(x):x\in T\}
$$
and
$$
R=\{x:x\in T,d_{H-S}(x)=0\}.
$$
We write $r=|R|$ and choose $x_1\in T$ such that
$d_{H-S}(x_1)=h_1$. If $T\setminus N_T[x_1]\neq\emptyset$, let
$$
h_2=\min\{d_{H-S}(x):x\in T\setminus N_T[x_1]\}.
$$
Thus, we have $0\leq h_1\leq h_2\leq a$ by the definition of $T$.

We shall consider various cases by the value of $r$ and derive a
contradiction in each case.

{\bf Case 1.} \ \ $r\geq2$.

Obviously, there exist $x,y\in R$ such that
$d_{H-S}(x)=d_{H-S}(y)=0$ and $xy\notin E(G)$. In terms of
$H=G-X$, Claim 5 and the condition of Theorem 1, we obtain
\begin{eqnarray*}
\frac{(a+b)n}{a+2b}&\leq&|N_G(x)\cup N_G(y)|\leq|N_H(x)\cup
N_H(y)|+|X|\\
&\leq&d_{H-S}(x)+d_{H-S}(y)+|S|+|X|=|S|+|X|<\frac{(a+b)n}{a+2b},
\end{eqnarray*}
which is a contradiction.

{\bf Case 2.} \ \ $r=1$.

Clearly, $h_1=0$ and $|N_T[x_1]|=1$. According to Claim 3, $r=1$
and $|N_T[x_1]|=1$, we have $T\setminus N_T[x_1]\neq\emptyset$ and
$1\leq h_2\leq a$. Choose $x_2\in T\setminus N_T[x_1]$ such that
$d_{H-S}(x_2)=h_2$. It is easy to see that $x_1x_2\notin E(G)$.
According to $H=G-X$ and the condition of Theorem 1, we have
\begin{eqnarray*}
\frac{(a+b)n}{a+2b}&\leq&|N_G(x_1)\cup N_G(x_2)|\leq|N_H(x_1)\cup
N_H(x_2)|+|X|\\
&\leq&d_{H-S}(x_1)+d_{H-S}(x_2)+|S|+|X|=h_2+|S|+|X|,
\end{eqnarray*}
which implies
$$
|S|\geq\frac{(a+b)n}{a+2b}-h_2-|X|.\eqno(2)
$$

Note that $|T\setminus N_T[x_1]|=|T|-1$. Combining this with
$|S|+|T|+|X|\leq n$, (2), Claim 1, $b\geq a\geq1$, $1\leq h_2\leq
a$ and
$n\geq\frac{(a+2b)(2a+2b-3)+1}{b}>\frac{(a+2b)(2a+2b-3)}{b}$, we
obtain
\begin{eqnarray*}
\delta_H(S,T)&=&b|S|+d_{H-S}(T)-a|T|=b|S|+d_{H-S}(N_T[x_1])+d_{H-S}(T\setminus N_T[x_1])-a|T|\\
&=&b|S|+d_{H-S}(T\setminus N_T[x_1])-a|T|\geq b|S|+h_2(|T|-1)-a|T|\\
&=&b|S|-(a-h_2)|T|-h_2\geq b|S|-(a-h_2)(n-|S|-|X|)-h_2\\
&=&(a+b-h_2)|S|-(a-h_2)n+(a-h_2)|X|-h_2\\
&\geq&(a+b-h_2)(\frac{(a+b)n}{a+2b}-h_2-|X|)-(a-h_2)n+(a-h_2)|X|-h_2\\
&=&(a+b-h_2)(\frac{(a+b)n}{a+2b}-h_2)-(a-h_2)n-b|X|-h_2\\
&\geq&(a+b-h_2)(\frac{(a+b)n}{a+2b}-h_2)-(a-h_2)n-\frac{b^{2}n}{a+2b}-h_2\\
&=&h_2^{2}+(\frac{bn}{a+2b}-a-b-1)h_2\\
&>&h_2^{2}+(\frac{(a+2b)(2a+2b-3)}{a+2b}-a-b-1)h_2\\
&=&h_2^{2}+(a+b-4)h_2\geq h_2^{2}-2h_2=(h_2-1)^{2}-1\geq-1,
\end{eqnarray*}
which contradicts (1).

{\bf Case 3.} \ \ $r=0$.

If $h_1=a$, then by (1) we obtain
$-1\geq\delta_H(S,T)=b|S|+d_{H-S}(T)-a|T|\geq
b|S|+h_1|T|-a|T|=b|S|\geq0$, which is a contradiction. Thus, we
have
$$
1\leq h_1\leq a-1.\eqno(3)
$$
We now prove the following claim.

{\bf Claim 6.} \ \ $T\setminus N_T[x_1]\neq\emptyset$.

{\bf Proof.} \ Suppose that $T=N_T[x_1]$. Then from (3) we have
$$
|T|=|N_T[x_1]|\leq|N_{H-S}[x_1]|=d_{H-S}(x_1)+1=h_1+1\leq a,
$$
which contradicts Claim 3.

In view of Claim 6, there exists $x_2\in T\setminus N_T[x_1]$ such
that $d_{H-S}(x_2)=h_2$. Obviously, $x_1x_2\notin E(G)$. According
to the condition of Theorem 1, we obtain
\begin{eqnarray*}
\frac{(a+b)n}{a+2b}&\leq&|N_G(x_1)\cup N_G(x_2)|\leq|N_H(x_1)\cup
N_H(x_2)|+|X|\\
&\leq&d_{H-S}(x_1)+d_{H-S}(x_2)+|S|+|X|=h_1+h_2+|S|+|X|,
\end{eqnarray*}
that is,
$$
|S|\geq\frac{(a+b)n}{a+2b}-h_1-h_2-|X|.\eqno(4)
$$
It is easy to see that
$$
|N_T[x_1]|\leq|N_{H-S}[x_1]|=d_{H-S}(x_1)+1=h_1+1.\eqno(5)
$$

Using $1\leq h_1\leq h_2\leq a$, $|S|+|T|+|X|\leq n$, (4), (5) and
Claim 1, we have
\begin{eqnarray*}
\delta_H(S,T)&=&b|S|+d_{H-S}(T)-a|T|=b|S|+d_{H-S}(N_T[x_1])+d_{H-S}(T\setminus N_T[x_1])-a|T|\\
&\geq&b|S|+h_1|N_T[x_1]|+h_2(|T|-|N_T[x_1]|)-a|T|=b|S|-(h_2-h_1)|N_T[x_1]|-(a-h_2)|T|\\
&\geq&b|S|-(h_2-h_1)(h_1+1)-(a-h_2)(n-|S|-|X|)\\
&=&(a+b-h_2)|S|-(h_2-h_1)(h_1+1)-(a-h_2)n+(a-h_2)|X|\\
&\geq&(a+b-h_2)(\frac{(a+b)n}{a+2b}-h_1-h_2-|X|)-(h_2-h_1)(h_1+1)-(a-h_2)n+(a-h_2)|X|\\
&=&(a+b-h_2)(\frac{(a+b)n}{a+2b}-h_1-h_2)-(h_2-h_1)(h_1+1)-(a-h_2)n-b|X|\\
&\geq&(a+b-h_2)(\frac{(a+b)n}{a+2b}-h_1-h_2)-(h_2-h_1)(h_1+1)-(a-h_2)n-\frac{b^{2}n}{a+2b}\\
&=&\frac{bn}{a+2b}h_2-(a+b-h_2)(h_1+h_2)-(h_2-h_1)(h_1+1),
\end{eqnarray*}
that is,
$$
\delta_H(S,T)\geq\frac{bn}{a+2b}h_2-(a+b-h_2)(h_1+h_2)-(h_2-h_1)(h_1+1).\eqno(6)
$$

Let
$F(h_1,h_2)=\frac{bn}{a+2b}h_2-(a+b-h_2)(h_1+h_2)-(h_2-h_1)(h_1+1)$.
Thus, by (3) we have
$$
\frac{\partial F(h_1,h_2)}{\partial
h_1}=-(a+b-h_2)-(-h_1-1+h_2-h_1)=-(a+b)+2h_1+1\leq-(a+b)+2(a-1)+1\leq-1.
$$
Combining this with $1\leq h_1\leq h_2\leq a$, we obtain
$$
F(h_1,h_2)\geq F(h_2,h_2).\eqno(7)
$$

In terms of (6), (7), $1\leq h_2\leq a$ and
$n\geq\frac{(a+2b)(2a+2b-3)+1}{b}>\frac{(a+2b)(2a+2b-3)}{b}$, we
have
\begin{eqnarray*}
\delta_H(S,T)&\geq&F(h_1,h_2)\geq
F(h_2,h_2)=\frac{bn}{a+2b}h_2-2(a+b-h_2)h_2\\
&>&\frac{(a+2b)(2a+2b-3)}{a+2b}h_2-2(a+b-h_2)h_2\\
&=&h_2(2h_2-3)\geq-1,
\end{eqnarray*}
which contradicts (1).

From all the cases above, we deduced the contradictions. Hence,
$H$ has a fractional $[a,b]$-factor, that is, $G$ is fractional
ID-$[a,b]$-factor-critical. The proof of Theorem 1 is complete.

\section{Remarks}
{\bf Remark 1.} \ In Theorem 1, the bound in the condition
$$
|N_G(x)\cup N_G(y)|\geq\frac{(a+b)n}{a+2b}
$$
is sharp. We can show this by constructing a graph
$G=(at)K_1\vee(bt)K_1\vee(bt+1)K_1$, where $t$ is sufficiently
large positive integer. It is easy to see that
$|V(G)|=n=(a+2b)t+1$ and
$$
\frac{(a+b)n}{a+2b}>|N_G(x)\cup
N_G(y)|=(a+b)t=(a+b)\cdot\frac{n-1}{a+2b}=\frac{(a+b)n}{a+2b}-\frac{a+b}{a+2b}>\frac{(a+b)n}{a+2b}-1
$$
for each pair of nonadjacent vertices $x,y$ of $(bt+1)K_1\subset
G$. Set $X=(bt)K_1$. Clearly, $X$ is an independent set of $G$.
Put $H=G-X=(at)K_1\vee(bt+1)K_1$, $S=(at)K_1$ and $T=(bt+1)K_1$.
Then $|S|=at$, $|T|=bt+1$ and $d_{H-S}(T)=0$. Thus, we have
\begin{eqnarray*}
\delta_H(S,T)&=&b|S|+d_{H-S}(T)-a|T|\\
&=&abt-a(bt+1)=-a<0.
\end{eqnarray*}
In terms of Lemma 2.1, $H$ has no fractional $[a,b]$-factor.
Hence, $G$ is not fractional ID-$[a,b]$-factor-critical.

\medskip

\medskip

{\bf Remark 2.} \ We show that the bound on minimum degree
$\delta(G)\geq\frac{bn}{a+2b}+a$ in Theorem 1 is also best
possible. Consider a graph $G$ constructed from $btK_1$,
$(at-1)K_1$, $\frac{bt}{2}K_2$ and $K_1$ as follows: let
$\{x_1,x_2,\cdots,x_{a-1}\}\subset (at-1)K_1$ and $K_1=\{u\}$,
where $t$ is sufficiently large positive integer and $bt$ is even.
Set $V(G)=V(btK_1\cup(at-1)K_1\cup\frac{bt}{2}K_2\cup\{u\})$ and
$E(G)=E(btK_1\vee(at-1)K_1\vee\frac{bt}{2}K_2)\cup
E(btK_1\vee\{u\})\cup\{ux_i:i=1,2,\cdots,a-1\}$. It is easily seen
that $|N_G(x)\cup N_G(y)\}|\geq\frac{(a+b)n}{a+2b}$ for each pair
of nonadjacent vertices $x,y$ of $G$, $n=(a+2b)t$ and
$\delta(G)=\frac{bn}{a+2b}+a-1$. Let $X=btK_1$. It is easy to see
that $X$ is an independent set of $G$. Set $H=G-X$, then
$\delta(H)=d_H(u)=a-1$. Clearly, $H$ has no fractional
$[a,b]$-factor, that is, $G$ is not fractional
ID-$[a,b]$-factor-critical.

\end{document}